\documentclass[10pt]{article}         % flush right qed marks, e.g. at end of proof
\usepackage{pgfplots}

\usepackage{geometry}
\usepackage{graphicx}
\usepackage{amssymb}
\usepackage{epstopdf}
\usepackage{amsmath,amsthm}
\usepackage{graphicx}
\usepackage{mathtools}
\usepackage{graphicx}
\usepackage{fix-cm}
\usepackage{pgf,tikz}
\usepackage{latexsym}
\usepackage{amsmath,amsfonts}
\usepackage{multirow}
\usepackage{latexsym}
\usepackage[english]{babel}
\usepackage{dcolumn}
\usepackage{color}
\usepackage{url}
\usepackage{algorithm2e}
\usepackage[justification=centering,labelfont=bf,font={small}]{caption}

\newcommand{\ignore}[1]{}

\newtheorem{theorem}{Theorem}
\newtheorem{corollary}{Corollary}

\newcommand{\oQ}{{\mathbf Q}}
\newcommand{\oR}{{\mathbb R}}

\newcommand{\oN}{{\mathbb N}}

\newcommand{\bK}{{\mathbf K}}

\begin{document}

\title{Distributionally robust optimization with polynomial densities: theory, models and algorithms} %theoretical properties and modeling approaches}
%\subtitle{Comparison of measure-based and simulated annealing bounds}

\author{Etienne de Klerk\thanks{Tilburg University \& TU Delft, The Netherlands; E.deKlerk@uvt.nl}    \and Daniel Kuhn\thanks{EPFL, Lausanne, Switzerland; daniel.kuhn@epfl.ch}  \and Krzysztof Postek\thanks{Erasmus University Rotterdam, The Netherlands; postek@ese.eur.nl}}

%\institute{Etienne de Klerk \at
%              Tilburg University \& TU Delft,
%               The Netherlands \\
%                             \email{E.deKlerk@uvt.nl} \and
%              Daniel Kuhn \at EPFL, Lausanne, Switzerland \\ \email{daniel.kuhn@epfl.ch} \and
%              Krzysztof Postek \at Erasmus University Rotterdam, The Netherlands \\ \email{postek@ese.eur.nl}
%              }

%\date{Received: date / Accepted: date}
% The correct dates will be entered by the editor

\maketitle

\begin{abstract}
{
In distributionally robust optimization the probability distribution of the uncertain problem parameters is itself uncertain, and a fictitious adversary, {\em e.g.}, nature, chooses the worst distribution from within a known ambiguity set. A common shortcoming of most existing distributionally robust optimization models is that their ambiguity sets contain pathological discrete distribution that give nature too much freedom to inflict damage. We thus introduce a new class of ambiguity sets that contain only distributions with sum-of-squares polynomial density functions of known degrees. We show that these ambiguity sets are highly expressive as they conveniently accommodate distributional information about higher-order moments, conditional probabilities, conditional moments or marginal distributions. Exploiting the theoretical properties of a measure-based hierarchy for polynomial optimization due to Lasserre [{\em SIAM J.~Optim.}~21(3) (2011), pp.~864--885], we prove that certain worst-case expectation constraints are computationally tractable under these new ambiguity sets.
We showcase the practical applicability of the proposed approach in the context of a stylized portfolio optimization problem and a risk aggregation problem of an insurance company.

}

{

% \PACS{PACS code1 \and PACS code2 }
%\subclass{90C22\and 90C26  \and 90C15}
}
\end{abstract}

\noindent
{\bf Keywords:} {distributionally robust optimization \and semidefinite programming \and sum-of-squares polynomials \and generalized eigenvalue problem}

\noindent
{\bf AMS classification:} 90C22, 90C26, 90C15

\section{Introduction}
\label{secintro}
Since George Dantzig's 1955 paper on linear programming under uncertainty \cite{Dantzig1955}, the field of stochastic programming has developed numerous methods for solving optimization problems that depend on uncertain parameters governed by a {\em known} probability distribution, see, {\em e.g.}, \cite{Birge1997,Prekopa1995,Shapiro2009}.
Stochastic programming usually aims to minimize a probability functional such as the expected value, a percentile or the conditional value-at-risk of a given cost function. In practice, however, the distribution that is needed to evaluate this probability functional is at best indirectly observable through independent training samples.
Thus, the stochastic programming approach is primarily useful when there is abundant training data. If data is scarce or absent, on the other hand, it may be more adequate to use a robust optimization approach, which models the uncertainty through the set of all possible (or sufficiently likely) uncertainty realizations and minimizes the worst-case costs. Robust optimization is the appropriate modeling paradigm for safety-critical applications with little tolerance for failure and has been popularized in the late 1990's, when it was discovered that robust optimization models often display better tractability properties than stochastic programming models~\cite{BenTal2009}. Distributionally robust optimization is a hybrid approach that attempts to salvage the tractability of robust optimization while maintaining the benefits of (limited) distributional information. In this context, uncertainty is modeled through an {\em ambiguity set}, that is, a family of typically infinitely many different distributions that are consistent with the available training data or any prior distributional information, and the objective is to minimize the worst-case expected costs across all distributions in the ambiguity set. A distributionally robust newsvendor model that admits an analytical solution has been investigated as early as in 1958 \cite{Scarf1958}, and the theoretical properties of distributionally robust linear programs were first studied in 1966 \cite{Zackova1966}.
Interest in distributionally robust optimization  has also been fuelled by important applications in finance \cite{Pflug2007,Pflug2012}.
However, only recently it was recognized that many distributionally robust optimization problems of practical relevance can actually be solved in polynomial time. Tractability results are available both for {\em moment ambiguity sets}, which contain all distributions that satisfy a finite number of moment conditions \cite{Delage2010,Goh2010,Wiesemann2014}, as well as for {\em metric-based ambiguity sets}, which contain all distributions within a prescribed distance from a nominal distribution with respect to some probability metric~\cite{BenTal2013,Esfahani2017}. In all these cases, the extremal distributions that determine the worst-case expectation are discrete, and the number of their discretization points is often surprisingly small, {\em e.g.}, proportional to the number of moment constraints. As these unnatural discrete distributions are almost always inconsistent with the available training samples, distributionally robust optimization models with moment and metric-based ambiguity sets are often perceived as overly pessimistic.

In an attempt to mitigate the over-conservatism of traditional distributionally robust optimization, several authors have studied moment ambiguity sets that require their member distributions to satisfy additional structural properties such as symmetry, unimodality, monotonicity or smoothness etc. By leveraging ideas from Choquet theory and polynomial optimization, it has been shown that the resulting distributionally robust optimization problems admit hierarchies of increasingly accurate semidefinite programming (SDP) bounds \cite{Popescu2005}. An exact SDP reformulation for the worst-case probability of a polytope with respect to all unimodal distributions with known first and second moments is derived in \cite{VanParys2016b}, while second-order conic reformulations of distributionally robust individual chance constraints with moment and unimodality information are reported in \cite{Li2017}. For a survey of recent results on distributionally robust uncertainty quantification and chance constrained programming problems with moment and structural information we refer to~\cite{Hanasusanto2015c}. Even though unimodality or monotonicity conditions eliminate all discrete distributions from a moment ambiguity set, the extremal distributions that critically determine all worst-case expectations remain pathological. For example, all extremal {\em unimodal} distributions are supported on line segments emanating from a single point in space (the mode) and thus fail to be absolutely continuous with respect to the Lebesgue measure. Thus, the existing distributionally robust optimization models with structural information remain overly conservative. This observation motivates us to investigate a new class of ambiguity sets that contain only distributions with non-degenerate polynomial density functions.

This paper aims to study \textit{worst-case expectation constraints} of the form \begin{equation} \label{eq.wcef} \sup\limits_{\mathbb{P} \in \mathcal{P}} \mathbb{E}_\mathbb{P} f(x,z) \leq 0, \end{equation} where $x\in \mathbb R^n$ is a decision vector, $z\in\mathbb R^m$ is an uncertain parameter governed by an ambiguous probability distribution $\mathbb{P} \in \mathcal{P}$, and $f(x,z)$ is an uncertainty-affected constraint function that can be interpreted as a cost. In words, the constraint~\eqref{eq.wcef} requires that the expected cost of the decision $x$ be non-positive for every distribution in the ambiguity set $\mathcal P$. Throughout the paper we will assume that $f(x,z)$ depends polynomially on $z$ and that each distribution $\mathbb P\in\mathcal P$ admits a sum-of-squares (hence non-negative) polynomial density function $h(z)$ with respect to some prescribed reference measure $\mu$ on $\mathbb R^m$ ({\em e.g.}, the Lebesgue measure). Imposing an upper bound on the polynomial degree of $h(z)$ thus yields a finite-dimensional parameterization of the ambiguity set $\mathcal P$.
Moreover, many popular distributional properties can be expressed through linear constraints on the coefficients of $h(z)$ and are thus conveniently accounted for in the definition of $\mathcal P$. Examples include moment bounds, probability bounds for certain subsets of $\mathbb R^m$, bounds on conditional tail probabilities and marginal distribution conditions. Note that by fixing the marginal distributions of all components of $z$, the worst-case expectation problem on the left-hand side of~\eqref{eq.wcef} reduces to a Fr\'echet problem that seeks the worst-case copula of the uncertain parameters.

By leveraging a measure-based hierarchy for polynomial optimization due to Lasserre~\cite{Las11}, we will demonstrate that the subordinate worst-case probability problem in~\eqref{eq.wcef} admits an exact SDP reformulation. Under mild additional conditions on $f(x,z)$, we will further prove that the feasible set of the constraint~\eqref{eq.wcef} admits a polynomial-time separation oracle. Moreover, we will analyze the convergence of the worst-case expectation in~\eqref{eq.wcef} as the polynomial degree of $h(z)$ tends to infinity, and we will illustrate the practical use of the proposed approach through numerical examples.

More succinctly, the main contributions of this paper can be summarized as follows:
\begin{itemize}
\item[\em (i)] \emph{Modeling power:} We introduce a new class of ambiguity sets containing distributions that admit sum-of-squares polynomial density functions of degree at most $2r$, $r\in\mathbb N$, with respect to a given reference measure. Ambiguity sets of this type are highly expressive as they conveniently accommodate distributional information about higher-order moments, conditional probabilities or conditional moments. They also allow the modeler to prescribe (not necessarily discrete) marginal distributions that must be matched exactly by all distributions in the ambiguity set.
\item[\em (ii)] \emph{Computational tractability:} We identify general conditions under which the worst-case expectations over the new ambiguity sets can be reformulated {\em exactly} as tractable SDPs with $\mathcal O{n+r \choose r}$ variables. We also propose an efficient heuristic for computing the worst-case expectations {\em approximately} by solving a sequence of significantly smaller SDPs. Finally, we delineate conditions under which the feasible sets of the worst-case expectation constraints admit a polynomial-time separation oracle and thus lend themselves to efficient optimization via the ellipsoid method.
\item[\em (iii)] \emph{Convergence analysis:} We demonstrate that, as $r$ tends to infinity, the worst-case expectations over the new ambiguity sets converge monotonically to classical worst-case expectations over larger ambiguity sets that relax the polynomial density requirement. At the same time, the extremal density functions converge to  pathological discrete worst-case distributions characteristic for classical moment ambiguity sets without restrictions on the density functions.
\item[\em (iv)] \emph{Numerical results:} We showcase the practical applicability of the proposed approach in the context of a stylyzed portfolio optimization problem and a simple Fr\'echet problem inspired by~\cite{VanParys2016} that models the risk aggregation problem of an insurance company.
\end{itemize}

The intimate relation between polynomial optimization and the problem of moments has already been exploited in several papers on distributionally robust optimization. For example, ideas from polynomial optimization give rise to SDP bounds on the probability of a semi-algebraic set \cite{Bertsimas2005} or the expected value of a piecewise polynomial \cite{Zuluaga2005} across all probability distributions satisfying a given set of moment constraints. These SDP bounds are tight in the univariate case or if only marginal moments are specified. Otherwise, one may obtain hierarchies of asymptotically tight SDP bounds. As an application, these techniques can be used to derive bounds on the prices of options with piecewise polynomial payoff functions, based solely on the knowledge of a few moments of the underlying asset prices \cite{Bertsimas2002}.
Moreover, asymptotically tight SDP bounds that account for both moment and structural information are proposed in \cite{Popescu2005}. However, all these approaches differ from our work in that the ambiguity sets have discrete or otherwise degenerate extremal distributions.

Distributionally robust polynomial optimization problems over non-degenerate polynomial density functions that are close to a nominal density estimate (obtained, {\em e.g.}, via a Legendre series density estimator) in terms of the $L_2$-distance are considered in \cite{Mevissen2013}. In this work the non-negativity of the candidate density functions is not enforced explicitly, which considerably simplifies the problem and may be justified if the distance to the nominal density is sufficiently small. It is shown that the emerging distributionally robust optimization problems are equivalent to deterministic polynomial optimization problems that are not significantly harder than the underlying nominal problem and can be addressed by solving a sequence of tractable SDP relaxations.

Distributionally robust chance constraints with ambiguity sets containing all possible mixtures of a given parametric distribution family are studied in \cite{Las11}. The mixtures are encoded through a probability density function on a compact parameter space. The authors propose an asymptotically tight SDP hierarchy of inner approximations for the feasible set of the distributionally robust chance constraint.
In contrast, we explicitly represent all probability distributions in the ambiguity set through polynomial density functions that can capture a wide range of distributional features.

The remainder of this paper develops as follows.
Section~\ref{sec:Lasserre hierarchy} reviews Lasserre's measure-based approach to polynomial optimization, which is central to this paper. A major drawback of the resulting SDP hierarchies is their limited scalability. This prompts us to devise an efficient heuristic solution algorithm in Section~\ref{section.heuristic}. Sections~\ref{sec:expected feasibility} and~\ref{sec:moment problem} develop SDP hierarchies for worst-case expectation constraints of the form~\eqref{eq.wcef} with and without moment information, respectively, and investigate the convergence of the underlying worst-case expectations as the degree of the polynomial density functions tends to infinity.
Section~\ref{sec:modeling-power} highlights the modeling power of the proposed approach, while Section~\ref{sec:results} reports on numerical results for a portfolio design problem as well as a risk aggregation problem of an insurance company. Conclusions are drawn in Section~\ref{sec:conclusions}.

\section{Lasserre's measure-based hierarchy for polynomial optimization}
\label{sec:Lasserre hierarchy}
In what follows, we denote by $x^\alpha :=\prod_{i=1}^nx_i^{\alpha_i}$ the monomial of the variables  $x = (x_1,\ldots,x_n)$ with respective exponents $\alpha = (\alpha_1,\ldots,\alpha_n) \in \oN_0^n$, and we define $N(n,r):=\{\alpha\in \oN_0^n: \sum_{i=1}^n\alpha_i\le r\}$ as the set of all exponents that give rise to monomials with degree at most~$r$.
We let $\Sigma[x]$ denote the set of all sum-of-squares (SOS) polynomials in the variables $x$, and we define $\Sigma[x]_r$ as the subset of all SOS polynomials with degree at most $2r$. %\smallskip\noindent

Now consider the polynomial global optimization problem
\begin{equation}
\label{eq:global-opt}
p_{\min,\bK} := \min_{x \in \bK} p(x) = \min_{x \in \bK} \sum_{\alpha \in N(n,d)} p_\alpha x^\alpha,
\end{equation}
where $p(x) = \sum_{\alpha \in N(n,d)} p_\alpha x^\alpha$ is an $n$-variate polynomial of degree $d$, and $\bK \subset \mathbb{R}^n$ a closed set with nonempty interior. (We assume existence of a global minimizer.)

We also assume that the moments of a finite Borel measure $\mu$ supported on $\bK$ are known in the sense that they are either available in closed-form or efficiently computable. To be clear,  we view a finite Borel measure $\mu$ on $\mathbb{R}^n$ as a nonnegative
set function defined on the Borel $\sigma$-algebra of $\mathbb{R}^n$. (Recall that the Borel $\sigma$-algebra is generated by all open sets in $\mathbb{R}^n$.)
By definition, $\mu$ must satisfy $\mu(\emptyset) = 0$ and $\mu(\cup_{i=1}^\infty S_i) = \sum_{i=1}^\infty \mu(S_i)$ for any collection of disjoint, measurable sets $S_i \subset \mathbb{R}^n$, $i\in\mathbb N$, and $\mu(\mathbb{R}^n) < \infty$.
The support of $\mu$, denoted by $\mbox{supp}(\mu)$, is defined as the smallest closed
set $\bK$ so that $\mu\left(\mathbb{R}^n\setminus \bK\right) = 0$.

 We denote the (known) moments of $\mu$ by
\begin{equation}\label{eq:moments}
m_{\alpha}(\mathbf{K}):=\int_{\mathbf{K}}x^{\alpha}d\mu(x)\ \ \ \text{ for } \alpha\in \oN_0^n.
\end{equation}
Lasserre \cite{Las11} introduced the following upper bound on $p_{\min,\bK}$,
\begin{eqnarray}
\underline{p}_{\mathbf{K}}^{(r)} &:=& \min_{h \in \Sigma_r} \left\{\int_{\mathbf{K}}p(x)h(x) d\mu(x) \; : \; \int_{\mathbf{K}}h(x) d\mu(x)=1 \right\} \label{eq:pminr} \\
&=& \min_{h \in \Sigma[x]_r} {\mathbb E}_{x \sim (\mathbf{K},h)} [p(x)], \nonumber
\end{eqnarray}
where $r$ is a fixed integer, and $x \sim (\mathbf{K},h)$ indicates that $x$ is a random vector supported on $\mathbf{K}$ that is governed by the probability measure $h\cdot d\mu$. It is known that if $\mu$ is the Lebesgue measure, then $\underline{p}_{\mathbf{K}}^{(r)}$ is equal to the the smallest generalized eigenvalue of the system
\begin{equation}
\label{eq:gen eig}
Av = \lambda Bv,
\end{equation}
with $v \neq 0$, where the symmetric matrices $A$ and $B$ are of size ${n + r \choose r}$ with rows and columns  indexed by $N(n,r)$,
and
\begin{equation}
\label{matrices A and B}
A_{\alpha, \beta} = \sum_{\delta \in N(n,d)} p_\delta m_{\alpha+\beta+\delta}(\mathbf{K}),
\quad B_{\alpha, \beta} = m_{\alpha+\beta}(\mathbf{K}) \quad \text{ for }  \alpha, \beta \in {N}(n,r).
\end{equation}
Lasserre \cite{Las11} establishes conditions on $\mu$ and $\bK$ so that $\lim_{r\rightarrow \infty} \underline{p}_{\mathbf{K}}^{(r)} = p_{\min,\bK}$, and the rate of
convergence was subsequently studied in \cite{KLS_MPA,KL_MOR_2017,KL_2018} for special choices of $\mu$ and $\bK$.
The most general condition under which convergence holds, as shown in \cite[Theorem 2.2]{Las12}, is when $\bK$ is closed with
nonempty interior, and the moments of $\mu$ on $\bK$ satisfy the following conditions:
\begin{equation}\label{eq:Carleman}
 \int_{\mathbf{K}}x_i^{2k}d\mu(x) \le (2k)!M \;\;\; \forall i\in \{1,\ldots,n\},k \in \mathbb{N},
\end{equation}
for some $M > 0$.
For example, if one defines $\mu$  in terms of a finite Borel measure $\varphi$ with supp$(\varphi) = \bK$ via
\begin{equation}\label{eq:unbounded measure}
  d\mu(x) = \exp\left(-|x_1|-\ldots -|x_n|\right)d\varphi(x),
\end{equation}
then this choice satisfies the conditions \eqref{eq:Carleman}; see \cite[\S3.2]{Las11}.

We summarize the known convergence results in Table~\ref{tab:convergence rates}.
\begin{table}[h!]
\begin{center}
\caption{Known rates of convergence for the Lasserre hierarchy \label{tab:convergence rates}}
\begin{tabular}{|c|c|c|c|}
  \hline
  $\bK \subset \mathbb{R}^n$  & $\underline{p}_{\mathbf{K}}^{(r)}- p_{\min,\bK}$ & measure $\mu$, supp$(\mu) = \bK$ & reference\\ \hline
  closed, nonempty interior & $o(1)$ & satisfies \eqref{eq:Carleman} & \cite{Las11}\\
  compact, nonempty interior & $o(1)$ & finite Borel measure & \cite{Las11}\\
  compact, satisfies interior cone condition &  $O\left(\frac{1}{\sqrt{r}}\right)$ & Lebesgue measure & \cite{KLS_MPA} \\
  convex body & $O\left(\frac{1}{r}\right)$ & Lebesgue measure & \cite{KL_MOR_2017}\\
  $[-1,1]^n$ & $\Theta\left(\frac{1}{r^2}\right)$ & $d\mu(x) = \prod_{i=1}^n (1-x_i^2)^{-1/2}dx$  &\cite{KL_2018}\\
   \hline
\end{tabular}
\end{center}
\end{table}

\subsection{Examples of known moments}
The moments \eqref{eq:moments} are available in closed-form, for example, if $\mu$ is the Lebesgue measure and $\bK$ is an ellipsoid or triangulated polytope;
 see, {\em e.g.}, \cite{Las11} and \cite{KLS_MPA}.
% Indeed, when $\mathbf{K}$ is the $n$-dimensional canonical simplex $\Delta_n=\{x\in\oR^n_+:\sum_{i=1}^n x_i \le 1\}$,
% the unit hypercube $\oQ_n=[0,1]^n$,  or the unit ball $B_1(0)=\{x\in \oR^n: \|x\|\le 1\}$,  there exist explicit formulas for the moments
% $m_{\alpha}(\mathbf{K})$ for the Lebesgue measure.
For the canonical simplex, $\Delta_n=\{x\in\oR^n_+:\sum_{i=1}^n x_i \le 1\}$, we have
\begin{equation}\label{mc0}
m_{\alpha}(\Delta_n)={\prod_{i=1}^n\alpha_i!\over (\sum_{i=1}^n \alpha_i +n)!},
\end{equation}
see, {\em e.g.}, \cite[Equation (2.4)]{LZ01} or \cite[Equation (2.2)]{GM78}. One may trivially verify that the moments for the hypercube $\oQ_n=[0,1]^n$ are given by %$m_{\alpha}(\oQ_n)$:
\begin{eqnarray*}\label{galphac}
m_{\alpha}(\oQ_n)=\int_{\oQ_n}x^{\alpha}dx=\prod_{i=1}^n \int_0^1x_i^{\alpha_i}dx_i=\prod_{i=1}^n \frac{1}{\alpha_i+1}.
\end{eqnarray*}
The moments for the unit Euclidean ball are given by
\begin{eqnarray}\label{malphaball}
m_{\alpha}(B_1(0)) = \left\{ \begin{array}{ll}
%\frac{\pi^{n/2}\prod_{i=1}^n\left(\alpha_i-1\right)!!}{\Gamma\left(1+{n+|\alpha|\over 2}\right)2^{|\alpha|/2}} =
{\pi^{(n-1)/2}2^{(n+1)/ 2} \prod_{i=1}^n\left(\alpha_i-1\right)!! \over (n+\sum_{i=1}^n \alpha_i)!!}
&    \textrm{\quad if $\alpha_i$ is even for all $i$,}\\
0 & \textrm{\quad otherwise,}
\end{array} \right.
\end{eqnarray}
where the double factorial of any integer $k$ is defined through
\begin{eqnarray*}
%\begin{displaymath}
k!! = \left\{ \begin{array}{l@{\quad}l}
k\cdot(k-2)\cdots 3\cdot1 & \textrm{if $k>0$ is odd,}\\
k\cdot(k-2)\cdots 4\cdot2 & \textrm{if $k>0$ is even,}\\
1 & \textrm{if $k=0$ or $k=-1$.}
\end{array} \right.
%\end{displaymath}
\end{eqnarray*}
When $\bK$ is an ellipsoid, one may obtain the moments  from \eqref{malphaball} by applying an affine transformation of variables.
Another tractable support set that will become relevant in Section~\ref{sec:portfolio} of this paper is the knapsack polytope, that is, the intersection of a hypercube and a half-space;
the moments for this and related polytopes are derived in \cite{Marichal2008}. Finally, in Section \ref{sec:Frechet} we will work with the nonnegative orthant $\bK = \mathbb{R}^n_+$. Since $\bK$ is unbounded in this case, we need to introduce a measure of the form \eqref{eq:unbounded measure}.
A suitable choice that corresponds to $d\varphi(x) =   2^n \exp\left(-\sum_{i=1}^n x_i\right)dx$ in \eqref{eq:unbounded measure}, is
\[
d\mu(x) = \exp\left(-\sum_{i=1}^n x_i\right)dx.
\]
This is the exponential measure associated with the orthogonal Laguerre polynomials. For more information, the reader is
referred to \cite[\S3.2]{Las11}. We will also use another choice of measure for $\bK = \mathbb{R}^n_+$ in Section \ref{sec:Frechet}, namely the lognormal measure,
\begin{equation}
	\label{eq:lognormal0}
	d\mu(x) = \prod_{i=1}^n \frac{1}{x_i v_i \sqrt{2\pi}} \exp\left( - \frac{(\ln(x_i) - \bar z_i)^2}{2 v_i^2} \right)dx_i,
\end{equation}
where $\bar z_i$ and $v_i$ represent prescribed location and scale parameters, $i = 1,\ldots,n$.
The moments of $\mu$ are  given by
\begin{equation}\label{eq:lognormal moments}
  m_{\alpha}(\mathbf K) =  \prod_{i=1}^n \exp(\alpha_i \bar z_i + (\alpha_i v_i)^2 / 2).
\end{equation}
One may readily verify that these moments do \emph{not} satisfy the bounds on the moments in \eqref{eq:Carleman}.
When using this measure we are therefore not guaranteed convergence of the Lasserre hierarchy.

We stress that, even though these examples of known moments are limited, they include  typical sets that are routinely used in (distributionally) robust optimization
to represent uncertainty sets or supports, most notably budget uncertainty sets and ellipsoids.

\section{An efficient, heuristic implementation of a Lasserre-type hierarchy} \label{section.heuristic}
The drawback of solving problem \eqref{eq:pminr} is that it involves operations with matrices of order ${n+r \choose r}$ for increasing values of $r$.
Thus one is limited to relatively small values of $n$ and $r$.

In this section we describe a weaker hierarchy of bounds that is similar in spirit to the hierarchy in \eqref{eq:pminr}, but where the
sizes of the corresponding generalized eigenvalue problems remain the same at each level of the hierarchy.
Conceptually, the idea is to use the optimal density function, say $h \in  \Sigma_r$ at level $r$, to approximate the optimal density function
at a higher level in the hierarchy.

To explain the idea, consider again the global optimization problem~\eqref{eq:global-opt}, which minimizes an $n$-variate polynomial $p$ over a convex body $\bK$. We assume that the moments
of a prescribed reference measure $\mu$ supported on $\bK$ are known, and we denote these moments by
\[
m^\mu_{\alpha}(\mathbf{K}):=\int_{\mathbf{K}}x^{\alpha}d\mu(x)\ \ \ \text{ for } \alpha\in \oN^n_0,
\]
where we add the superscript $\mu$ to make the dependence on the reference measure explicit.
Next we compute the  upper bound $\underline{p}_{\mathbf{K}}^{(r)}$ as in \eqref{eq:pminr}, where $r$ is a fixed (small) integer. Denoting the resulting optimal density by $h(x) = \sum_{\beta \in N(n,2r)} h_\beta x^\beta \in \Sigma[x]_r$, we can then define a new probability measure~$\mu'$ on $\bK$ through
\begin{equation}
\label{eq:update density}
d\mu'(x) = h(x)\cdot d\mu(x).
\end{equation}
Note that we may obtain the moments of $\mu'$ from the moments of $\mu$ via
\begin{eqnarray*}
\label{eq:update moments}
m^{\mu'}_{\alpha}(\mathbf{K}) &=& \int_{\mathbf{K}}x^{\alpha}d\mu'(x)\\
&=& \int_{\mathbf{K}}x^{\alpha}h(x)d\mu(x)\\
&=& \sum_{\beta \in N(n,2r)} h_\beta\int_{\mathbf{K}}x^{\alpha+\beta}d\mu(x)\\
&=& \sum_{\beta \in N(n,2r)} h_\beta m^\mu_{\alpha+\beta}(\mathbf{K})\ \ \ \text{ for } \alpha\in \oN^n_0.\\
\end{eqnarray*}
Finally, one may now replace $\mu$ by $\mu'$ and repeat the same process $R$ times for some fixed $R\in\mathbb N$. The complete procedure is summarized in Algorithm \ref{alg:heuristic}.

\begin{center}
\begin{algorithm}[h!]
 \KwData{Polynomial $p$ of degree $d$, allowed degree $r\ge d$, integer order $R$; moments
 of some measure $\mu$ up to order $2r+d + 2(R-1)r = 2rR+d$, {\em i.e.}, the values
  $m^\mu_{\alpha}(\mathbf{K})$ for all $\alpha \in N(n,2rR+d)$  }
 \KwResult{Upper bound of order $R$ on $p_{\min,\bK}$}
 \For{$k\leftarrow 1$ \KwTo $R$}{
  Form the matrices $A$ and $B$ defined in \eqref{matrices A and B} \;
 Solve the generalized eigenvalue problem for $A$ and $B$ in \eqref{eq:gen eig} to obtain the optimal density $h \in \Sigma_r$ \;
 Define the measure $\mu'$ via \eqref{eq:update density} \;
 Obtain the moments of $\mu'$ via $m^{\mu'}_{\alpha}(\mathbf{K}) = \sum_{\beta \in N(n,2r)} h_\beta m^\mu_{\alpha+\beta}(\mathbf{K})$ for
 all $\alpha \in N(n,2r(R-k) + d)$ \;
 Replace $\mu \leftarrow \mu'$ and $m^{\mu}_{\alpha}(\mathbf{K}) \leftarrow m^{\mu'}_{\alpha}(\mathbf{K})$ for
 all $\alpha \in N(n,2r(R-k) + d)$ \;
 }
 \caption{Algorithm to compute the upper bound $\underline{p}_{\mathbf{K}}^{(r,R)}$ of order $R$ on $\underline{p}_{\mathbf{K}}^{(r)}$.} \label{alg:heuristic}
\end{algorithm}
\end{center}
The following remarks on this heuristic procedure are in order:
\begin{enumerate}
  \item
  The bound computed by the algorithm is no better than the bound $\underline{p}_{\mathbf{K}}^{(r\cdot R)}$,
  but is much cheaper to compute because, in each iteration, it only involves generalized eigenvalue problems of order ${n+r \choose r}$ for a small fixed integer $r$, {\em e.g.}, $r=4$.
    \item
  The bounds generated by the algorithm, as indexed by $R$, are not guaranteed to converge to $p_{\min,\bK}$ as
  $R\rightarrow \infty$. However, one may easily obtain a convergent variant (at a computational cost) by increasing the value $r$
  inside a given iteration, if no improvement in the upper bound is obtained in that iteration.
  \item
  One has to store and update a moment table indexed by $\alpha \in N(n,2r(R-k) + d)$ in  iteration~$k$, and the updating process involves simple linear algebra.
\end{enumerate}

We now give some numerical examples to indicate how this algorithm performs. We will consider the test functions
in Table \ref{tab:test}.

\begin{table}[h!]
\caption{Test functions, all with $n=2$, domain $\bK = [-1,1]^2$, and minimum $p_{\min,\bK} = 0$. \label{tab:test}}
\begin{center}
\begin{tabular}{|c|c|}\hline
Name & $p(x)$  \\ \hline
%Booth function& $(10x_1+20x_2-7)^2 + (20x_1+10x_2-5)^2$ & $2594$ & $2$ & yes \\ \hline
Matyas function& $26(x_1^2+x_2^2)-48x_1x_2$   \\ \hline
Motzkin polynomial& $64(x_1^4x_2^2+x_1^2x_2^4) - 48x_1^2x_2^2 +1$   \\ \hline
%Three-Hump Camel function& $\frac{5^6}{6}x_1^6-5^4\cdot 1.05x_1^4+50x_1^2+25x_1x_2+25x_2^2$ & $2048$ & $6$ &  no \\ \hline
\end{tabular}
 \end{center}
\end{table}

For ease of reference, we will denote the bound generated by Algorithm \ref{alg:heuristic} by
$\underline{p}_{\mathbf{K}}^{(r,R)}$. As this bound corresponds to a density function
of degree $rR$ with respect to the initial reference measure, it is natural to compare it to the
stronger, but more expensive, bound $\underline{p}_{\mathbf{K}}^{(r\cdot R)}$. The following table  lists different bounds for the case
$r\cdot R=20$, each corresponding to a density function of degree $40$.

\begin{table}[h!]
\begin{center}
\begin{tabular}{|c|c|c|c|c|c|}\hline
Function & $\underline{p}_{\mathbf{K}}^{(20)}$ & $\underline{p}_{\mathbf{K}}^{(10,2)}$ & $\underline{p}_{\mathbf{K}}^{(5,4)}$ & $\underline{p}_{\mathbf{K}}^{(4,5)}$ & $\underline{p}_{\mathbf{K}}^{(2,10)}$ \\ \hline
Matyas & 0.4811                                & 0.4989                   &            0.7285  & 0.9604 & 1.1070\\
Motzkin & 0.1817                               & 0.1969                                & 0.2907  &0.3603 & 0.4588 \\ \hline
 \end{tabular}
  \end{center}
  \end{table}
Note that the $\underline{p}_{\mathbf{K}}^{(r,R)}$ bounds with largest $r$ are the strongest in the examples, as one may expect.

\section{Distributionally robust constraints involving polynomial uncertainty}
\label{sec:expected feasibility}
We now consider a worst-case feasibility constraint of the form \eqref{eq.wcef}, where $z \in \mathbb{R}^m$
represents a random vector with a support $\bK \subset \mathbb{R}^m$, assumed to be closed and with nonempty interior.
 Assume that the constraint function $f(x,z)$ displays a polynomial dependence on $z$.
 In particular, assume that $f(x,z)=\sum_{\beta\in N(m,d)}f_{\beta}(x)z^\beta$ has degree
   $d$ in $z$, where the $f_\beta:\mathbb{R}^n \rightarrow \mathbb{R}$ are  functions of $x$ only.

If the ambiguity set $\mathcal P$ contains all distributions that have an SOS polynomial
density of degree at most $2r$, $r>1$, with respect to a a fixed, finite Borel measure $\mu$ supported on $\mathbf{K}$,
 then the worst-case expected feasibility constraint \eqref{eq.wcef} reduces to
\begin{eqnarray}\label{fminkreform2}
f_{\bK}^{(r)}(x) := \sup_{h\in\Sigma[z]_r}\left\{\int_{\mathbf{K}}f(x,z)h(z)d\mu(z) \ \ : \; \int_{\mathbf{K}}h(z)d\mu(z)=1\right\} \le 0.
\end{eqnarray}
 Formally speaking, we consider an ambiguity set of the form
\begin{equation}
\label{q:ambiguity-set}
\mathcal{P} = \left\{h\cdot d\mu \;:\; h\in\Sigma[z]_r, \; \int_{\mathbf{K}} h(z)d\mu(z)=1 \right\}.
\end{equation}
We assume that the moments of the measure $\mu$ on $\bK$ are available, and we again use the notation
\begin{equation*}\label{mack}
m_{\alpha}(\mathbf{K}):=\int_{\mathbf{K}}z^{\alpha}d\mu(z)\ \ \ \text{ for } \alpha\in \oN^m_0.
\end{equation*}
Expressing $h\in\Sigma[z]_{r}$ as $h(z)=\sum_{\alpha\in N(m,2r)}h_{\alpha}z^{\alpha}$,  the
left-hand-side of \eqref{fminkreform2} may be re-written~as
\begin{equation}\label{eqSDP}
\begin{array}{cl}
\displaystyle \sup_{h_\alpha : \alpha\in N(m,2r)} & \displaystyle \sum_{\beta\in N(m,d)}f_{\beta}(x)\sum_{\alpha\in N(m,2r)}h_{\alpha}m_{\alpha+\beta}(\mathbf{K})\\[1ex]
\text{ s.t. } & \displaystyle \sum_{\alpha\in N(m,2r)}h_{\alpha}m_{\alpha}(\mathbf{K})=1,\\[1ex]
& \displaystyle \sum_{\alpha\in N(m,2r)}h_{\alpha}z^{\alpha}\in\Sigma[z]_r.
\end{array}
\end{equation}
Since the sum-of-squares condition on $h$ is equivalent to a
linear matrix inequality in its coefficients~$h_\alpha$, problem~\eqref{eqSDP}
 constitutes a tractable semidefinite program (SDP) in $h_\alpha$, $\alpha\in N(m,2r)$, if $x$ is fixed.
  The next theorem establishes that we can also efficiently optimize over the feasible set of the constraint~\eqref{fminkreform2}
   whenever the coefficient functions $f_\beta$ are convex and $\bK \subset \mathbb{R}^m_+$.

\begin{theorem}
\label{thm:tractability}
Consider the constraint (\ref{fminkreform2}) and assume that all  $f_{\beta}$ are convex functions of $x$ whose subgradients are efficiently computable. Moreover, assume that $\bK \subset \mathbb{R}^m_+$. Then, the set of $x \in \mathbb{R}^n$ that satisfy (\ref{fminkreform2}) is convex and admits a polynomial-time separating hyperplane oracle.
\end{theorem}
\proof
We have to show that the function $f_{\bK}^{(r)}(x)$ from \eqref{fminkreform2} is convex in $x$.
We may rewrite the function as
\[
f_{\bK}^{(r)}(x) = \sup_{h \in \Sigma[z]_r}  \sum_{\beta\in N(m,d)} \mathbb{E}_{z \sim (\bK,h)}\left[ z^\beta\right]f_{\beta}(x)
\]
For each $h \in \Sigma[z]_r$, the function $\mathbb{E}_{z \sim (\bK,h)}\left[ z^\beta\right]f_{\beta}(x)$ is convex in $x$, since $\bK \subset \mathbb{R}^m_+$ implies
$\mathbb{E}_{z \sim (\bK,h)}\left[ z^\beta\right] \ge 0$.
Thus $f_{\bK}^{(r)}(x)$ is the point-wise supremum of an infinite collection of convex functions, and therefore convex itself
(see, {\em e.g.}, \cite[Theorem 5.5]{Rockafellar}).
Thus the set $\mathcal{C} := \{x \in \mathbb{R}^n \; | \; f_{\bK}^{(r)}(x) \le 0\}$ is convex.

If $\bar x \notin \mathcal{C}$, {\em i.e.}, $f_{\bK}^{(r)}(\bar x) > 0$, then we may construct a hyperplane that
separates $\bar x$  from $\mathcal{C}$ as follows. Let $\bar h \in \Sigma[z]_r$ be such that
\[
f_{\bar h}(\bar x) := \sum_{\beta\in N(m,d)} \mathbb{E}_{z \sim (\bK,\bar h)}\left[ z^\beta\right]f_{\beta}(\bar x) > 0.
\]
One may obtain such an $\bar h$ in polynomial time by solving the SDP \eqref{eqSDP} with fixed $x= \bar x$.
Now let $\partial f_{\bar h}(\bar x)$ denote a subgradient of $ f_{\bar h}$ at $\bar x$. (By assumption such a subgradient is available in polynomial time.)
By the definition of a subgradient, we now have
\[
\partial f_{\bar h}(\bar x)^T (x - \bar x) \le f_{\bar h}(x) - f_{\bar h}(\bar x) \le -f_{\bar h}(\bar x) \quad \forall x \in \mathcal{C}.
\]
The outer linear inequality now separates $\bar x$ from $\mathcal{C}$.
 \qed

Theorem~\ref{thm:tractability} implies that if all coefficient functions $f_{\beta}$ are convex,  one may optimize a convex function of $x$ over a set given by constraints of the type
(\ref{fminkreform2}) in polynomial time, {\em e.g.}, by using the ellipsoid method, provided that an initial ellipsoid is known that contains an
 optimal solution \cite{GroetschelLovaszSchrijver1988a}.

Finally, we point out that, due to the convergence properties of the Lasserre hierarchy, one recovers the usual robust counterpart (robust against the single worst-case realization of $z$ as in \cite{BenTal2009}) in the limit as $r$ tends to infinity.

\begin{theorem}
\label{thm:convergence}
Assume that $\bK\subset \mathbb{R}^n$ is closed with nonempty interior. Then, in the limit as $r \rightarrow \infty$, the constraint (\ref{fminkreform2}) reduces to the usual robust counterpart constraint
\[
\max_{z \in \bK} f(x,z) \le 0.
\]
More precisely, if $x \in \bK$ is fixed, and $(\bK,\mu)$ satisfies one of the assumptions in Table \ref{tab:convergence rates}, one has
\[
\lim_{r \rightarrow \infty} f_{\bK}^{(r)}(x) = \max_{z \in \bK} f(x,z).
\]
Moreover the rate of convergence is as given in Table \ref{tab:convergence rates}, depending on the choice of $(\bK,\mu)$.
\end{theorem}
\proof
For fixed $x$, the computation of $f_{\bK}^{(r)}(x)$ is an SDP problem of the form \eqref{eq:pminr}, and the required convergence result
therefore follows from the convergence of the Lasserre hierarchy \eqref{eq:pminr},  as summarized in Table \ref{tab:convergence rates}. \qed

\section{Approximate solution of the general problem of moments}
\label{sec:moment problem}
In applications it is often possible to inject moment information into the ambiguity set $\mathcal P$.
For example, if there is prior information about the location or the dispersion of the random vector~$z$,
one can include constraints on its mean vector or its covariance matrix into the definition of the ambiguity set.
Specifically, if it is known that ${\mathbb E}_{z \sim (\mathbf{K},h)} [z^{\beta_i}]=\gamma_i$ for some $\beta_i\in \mathbb N_0^n$ and $\gamma_i\in\mathbb R$ for $i=1,\ldots,p$, one can
restrict the ambiguity set~\eqref{q:ambiguity-set} by including the moment constraints
\begin{equation*}
	\int_{\mathbf{K}} z^{\beta_i} h(z) d\mu(z)=\sum_{\alpha\in N(m,2r)}h_{\alpha}m_{\alpha+\beta_i}(\mathbf{K})=\gamma_i \quad \forall i = 1,\ldots,p,
\end{equation*}
which reduce to simple linear equations for the coefficients $h_{\alpha}$, $\alpha\in N(m,2r)$, of the density function~$h$. In this setup, the maximization over the ambiguity set corresponds to a general problem of moments, see, {\em e.g.}, \cite{Shapiro2001}. We show how our approach may be used to solve this problem, and we will illustrate this through concrete examples in Section~\ref{sec:results}.

Throughout this section we assume that $\bK \subset \mathbb{R}^k$ is a nonempty closed set, while $f_0, f_1, \ldots, f_p$ are real-valued Borel-measurable functions on $\bK$. Moreover, we assume that $\mu$ is a finite Borel measure
on $\bK$ such that $f_0,\ldots,f_p$ are $\mu$-integrable.

\begin{theorem}
Let $\mathbf K_i$, $i=0,\ldots,p$, be Borel-measureable subsets of~$\bK$.
Then there exists an atomic Borel measure $\mu'$ on $\bK$ with a finite support of at most $p+2$ points so that
\[
\int_{\bK_i} f_i(z)d\mu(z) = \int_{\bK_i} f_i(z)d\mu'(z) \quad \forall i = 0,\ldots,p.
\]
\end{theorem}
\proof
This result is due to Rogosinsky \cite{rogosinski}, but an elementary proof  is given by Shapiro \cite[Lemma~3.1]{Shapiro2001}; see also Lasserre \cite{Las2008}.\qed

As a consequence one has the following result for the problem of moments.

\begin{corollary}
\label{cor:atomic moments}
Consider the problem of moments
\begin{equation}
\label{eq:prob moments}
	val:=\inf_{\mu \in \mathcal{P}} \left\{ \int_{\bK_0} f_0(z)d\mu(z) \; : \; \int_{\bK_i} f_i(z)d\mu(z) = b_i \;\forall i = 1,\ldots,p\right\},
\end{equation}
where $\mathcal{P}$ is the set of all Borel probability measures supported on $\bK$, and $\bK_i$ is a Borel-measurable subset of $\bK$ for each $i=0,\ldots,p$.
If the problem has a solution, it has a solution that is an atomic measure supported on at most $p+2$ points in $\bK$, i.e., a convex combination
of at most $p+2$ Dirac delta measures supported in $\bK$.
\end{corollary}

In what follows we show how the atomic measure solution, whose existence is guaranteed by Corollary~\ref{cor:atomic moments}, may be approximated arbitrarily well by SOS polynomial density functions.

\begin{theorem}
\label{th:moments sos}
  Consider problem
  \eqref{eq:prob moments} with the additional assumptions that $\bK \subset \mathbb{R}^n$ has nonempty interior
  and that the functions $f_0$, $f_1$, $\ldots$, $f_p$ are polynomials. Also assume that $\bK$ and $\mu$ satisfy one of the assumptions in Table \ref{tab:convergence rates}.
  Then, for any $\epsilon > 0$ there exists a $d \in \mathbb{N}$ and a probability density $h \in \Sigma_d[z]$ such that
  one has
  \begin{eqnarray*}
  % \nonumber % Remove numbering (before each equation)
    \int_{\bK_0} f_0(z)h(z)d\mu(z) &\in& (val-\epsilon,val+\epsilon) \\
      \int_{\bK_i} f_i(z)h(z)d\mu(z) &=& (b_i-\epsilon,b_i+\epsilon) \quad \forall i = 1,\ldots,p.
  \end{eqnarray*}
 Moreover, for the choices of $\bK$ and $\mu$ in Table \ref{tab:convergence rates} where a rate of convergence is known,
 one may bound $d$ in terms of $\epsilon$.
 For example, if $\bK$ is a convex body and $\mu$ the Lebesgue measure, then one may assume that $d = O(1/\epsilon^2)$.
\end{theorem}
\proof
Fix $a \in \bK$, and consider the polynomials $z \mapsto (f_i(z) - f_i(a))^2$, $i =0,\ldots,p$.
Moreover, let $p$ be a polynomial with global minimizer $a$ such that $p(a) = 0$ and
$p$ upper bounds all these polynomials  on $\bK$, {\em i.e.},
\begin{equation}\label{eq:dominated pols}
  p(z) \ge (f_i(z) - f_i(a))^2 \mbox{ for all } z \in \bK \mbox{ and all } i \in \{0,\ldots,p\}.
\end{equation}
For a given probability density $h \in \Sigma[z]_r$ with $\int_{\bK} h(z)d\mu(z) =1$, we denote as before
\[
\mathbb{E}_{z \sim (\bK,h)} [p(z)] = \int_{\bK} p(z) h(z)d\mu(z).
\]
Note that  by \eqref{eq:dominated pols} we have $\mathbb{E}_{z \sim (\bK,h)} [p(z)] \ge \mathbb{E}_{z \sim (\bK,h)} [(f_i(z) - f_i(a))^2]$.
Combining with Jensen's inequality, we therefore conclude
\begin{eqnarray*}\label{eq:Jensen}
 \left(\mathbb{E}_{z \sim (\bK,h)} [(f_i(z) - f_i(a))]\right)^2   &\le &  \mathbb{E}_{z \sim (\bK,h)} [(f_i(z) - f_i(a))^2]\\
   &\le & \mathbb{E}_{z \sim (\bK,h)} [p(z)].
\end{eqnarray*}
Recalling the notation of the Lasserre hierarchy from \eqref{eq:pminr}, we denote
$\underline{p}_{\mathbf{K}}^{(r)} =\min_{h \in \Sigma[z]_r} {\mathbb E}_{z \sim (\mathbf{K},h)} [p(z)]$.
If $\mu$ and $\bK$ satisfy one of the conditions from Table~\ref{tab:convergence rates}, one has
$\lim_{r \rightarrow \infty} \underline{p}_{\mathbf{K}}^{(r)} = 0$, with the rate of convergence as indicated in the table.
Thus, for any $\epsilon > 0$ there is a sufficiently large  $d \in \mathbb{N}$ such that
\[
\min_{h \in \Sigma[z]_r}  \left(\mathbb{E}_{z \sim (\bK,h)} [(f_i(z) - f_i(a))]\right)^2 \le {\epsilon} \quad \forall r \ge d, \; i \in \{0,\ldots,p\}.
\]
Letting $h^*$ denote the minimizer, one has
\[
\left|   \mathbb{E}_{z \sim (\bK,h^*)} [(f_i(z)] - f_i(a) \right|\le \sqrt{\epsilon}  \quad \forall r \ge d, \; i \in \{0,\ldots,p\}.
\]
To complete the proof, we simply have to associate $a$ with an atom of the optimal atomic distribution from Corollary \ref{cor:atomic moments}. \qed

As a consequence  of Theorem \ref{th:moments sos}, we may obtain approximate solutions to the generalized problem of moments
\eqref{eq:prob moments} by solving SDPs of the form:
\[
\min_{h \in \Sigma[z]_r} \left\{ \int_{\bK_0} f_0(z)h(z)d\mu(z) \; : \; \int_{\bK_i} f_i(z)h(z)d\mu(z) = [b_i-\epsilon,b_i+\epsilon] \;\forall i = 1,\ldots,p\right\},
\]
for given $r \in \mathbb{N}$ and $\epsilon \ge 0$, and we will do precisely that in the example of Section~\ref{sec:portfolio}.
We remark that these SDPs are different from the ones studied by Lasserre \cite{Las2008}, where an \emph{outer} approximation of the cone of finite Borel measures supported on $\bK$ is used, whereas we use an inner approximation.

\section{Modeling power}
\label{sec:modeling-power}
The ambiguity set $\mathcal P$ defined in~\eqref{q:ambiguity-set} contains all distributions supported on a convex body $\mathbf K$ that have an SOS polynomial density $h\in\Sigma[z]_r$ with respect to a prescribed reference measure $\mu$. For any fixed $x$, the worst-case expectation $f_{\bK}^{(r)}(x)$ on the left-hand-side of the worst-case feasibility constraint~\eqref{fminkreform2} can be computed efficiently by solving the SDP~\eqref{eqSDP}. The ambiguity set $\mathcal P$ admits several generalizations that preserve the SDP-representability of the worst-case expectation.

\paragraph{Moment information} As already discussed in Section~\ref{sec:moment problem}, conditions on (mixed)
 moment values of different random variables give rise to simple linear conditions on the polynomial coefficients of $h$.

\paragraph{Confidence information} If the random vector $z$ is known to materialize inside a given Borel set $\mathbf C\subset \mathbb R^m$ with probability $\gamma\in [0,1]$, we can add the condition ${\mathbb P}_{z \sim (\mathbf{K},h)} [z\in\mathbf C]=\gamma$ to the definition of the ambiguity set $\mathcal P$. Moreover, if the moments $m_\alpha(\mathbf K\cap \mathbf C)$ of the reference measure $\mu$ over $\mathbf K\cap \mathbf C$ are either available analytically or efficiently computable for all $\alpha \in N(m,2r)$, then this condition can be re-expressed as the following simple linear equation in the polynomial coefficients of $h$.
\begin{align*}
\int_{\mathbf{K}} \mathbf 1_{z\in\mathbf C}\, h(z) d\mu(z)= \sum_{\alpha\in N(m,2r)} h_{\alpha}m_{\alpha}(\mathbf{K} \cap \mathbf C) = \gamma
\end{align*}
Upper and lower bounds on ${\mathbb P}_{z \sim (\mathbf{K},h)} [z\in\mathbf C]$ can be handled similarly in the obvious manner. In the context of purely moment-based ambiguity sets, such probability bounds have been studied in~\cite{Wiesemann2014}.

\paragraph{Conditional probabilities} Given any two Borel sets $\mathbf C_1,\mathbf C_2\subset \mathbb R^m$ and a probability $\gamma\in [0,1]$, we can also enforce the condition ${\mathbb P}_{z \sim (\mathbf{K},h)} [z\in\mathbf C_2|z\in\mathbf C_1]=\gamma$ in the definition of $\mathcal P$. If the moments $m_\alpha(\mathbf K\cap \mathbf C_1)$ and $m_\alpha(\mathbf K\cap \mathbf C_1\cap \mathbf C_2)$ of the reference measure $\mu$ are either available analytically or efficiently computable for all $\alpha \in N(m,2r)$, then this condition can be re-expressed as
\begin{align*}
& \int_{\mathbf{K}} \mathbf 1_{z\in\mathbf C_12\cap \mathbf C_2}\, h(z) d\mu(z)= \gamma  \int_{\mathbf{K}} \mathbf 1_{z\in\mathbf C_1}\, h(z) d\mu(z) \\
\iff\quad & \sum_{\alpha\in N(m,2r)} h_{\alpha}\left( m_{\alpha}(\mathbf{K} \cap \mathbf C_1\cap \mathbf C_2 ) - \gamma\, m_{\alpha}(\mathbf{K} \cap \mathbf C_1 )\right)=0,
\end{align*}
which is again linear in the coefficients of $h$. Upper and lower bounds on conditional probabilities can be handled similarly.

\paragraph{Conditional moment information}  If it is known that ${\mathbb E}_{z \sim (\mathbf{K},h)} [z^{\beta}|\mathbf C]=\gamma$ for some $\beta\in \mathbb N_0^n$, Borel set $\mathbf C\subset \mathbb R^m$ and $\gamma\in\mathbb R$, while the moments $m_{\alpha + \beta}(\mathbf K\cap \mathbf C)$ of the reference measure $\mu$ over set $\mathbf K\cap \mathbf C$ are either available analytically or efficiently computable for all $\alpha \in N(m,2r)$, then one can add the following condition to the ambiguity set~$\mathcal P$, which is linear in the coefficients of $h$.
\begin{align*}
&	\int_{\mathbf{K}} z^{\beta} \mathbf 1_{z\in\mathbf C}\, h(z) d\mu(z)= \gamma \int_{\mathbf{K}} \mathbf 1_{z\in\mathbf C}\, h(z) d\mu(z) \\
\iff\quad & \sum_{\alpha\in N(m,2r)} h_{\alpha}\left( m_{\alpha+\beta}(\mathbf{K}\cap \mathbf C)  - \gamma\, m_{\alpha}(\mathbf{K} \cap \mathbf C )\right)=0,
\end{align*}

\paragraph{Multiple reference measures} The distributions in the ambiguity set $\mathcal P$ defined in~\eqref{q:ambiguity-set} depend both on the reference measure $\mu$ as well as the density function $h$. A richer ambiguity set can be constructed by specifying multiple reference measures $\mu_i$ with corresponding density functions $h^i\in\Sigma[z]_r$, $i=1,\ldots,p$. The distributions in the resulting ambiguity set are of the form $\sum_{i=1}^ph^i\cdot d\mu_i$. If the moments $m^i_\alpha(\mathbf K)$ of the reference measure $\mu_i$ over $\mathbf K$ are either available analytically or efficiently computable for all $\alpha \in N(m,2r)$ and $i=1,\ldots,p$, then the normalization constraint can be recast~as
\[
	\sum_{\alpha\in N(m,2r)} h^i_{\alpha}m^i_{\alpha}(\mathbf{K})=\gamma_i~\forall i=1,\ldots,p\quad \text{and} \quad \sum_{i=1}^p\gamma_i=1,
\]
where $\gamma=(\gamma_1,\ldots,\gamma_p)\ge 0$ constitutes an auxiliary decision vector. The resulting ambiguity set can be interpreted as a convex combination of $p$ ambiguity sets of the form~\eqref{q:ambiguity-set} and thus lends itself for modeling multimodality information; see, {\em e.g.}, \cite{Hanasusanto2015b}. In this case, $\gamma_i$ captures the probability of the $i$-th mode, which may itself be uncertain. Thus, $\gamma$ should range over a subset of the probability simplex, {\em e.g.}, a $\phi$-divergence uncertainty set of the type studied in~\cite{BenTal2013}.

\paragraph{Marginal distributions} It is often easier to estimate the marginal distributions of all $m$ components of a random vector $z$ instead of the full joint distribution. Marginal distribution information can also be conveniently encoded in ambiguity sets of the type~\eqref{q:ambiguity-set}. To see this, assume that the marginal distribution of $z_i$ is given by $\mu_i$ and is supported on a compact interval $\mathbf K_i\subset \mathbb R$, $i=1,\ldots,m$. In this case it makes sense to set $\mathbf K=\bigtimes_{i=1}^m \mathbf K_i$ and to define the reference measure $\mu$ through $d\mu=\prod_{i=1}^m d\mu_i$. Thus, $\mu$ coincides with the product of the known marginals. The requirement
\[
 \int_{\mathop{\bigtimes}_{j\neq i} \mathbf K_j} h(z) \prod_{j\neq i} d\mu_j(z_j) = 1 \quad \forall z_i\in\mathbf K_i,~\forall i = 1,\ldots,m
 \]
then ensures that the marginal distribution of $z_i$ under $h\cdot d\mu$ exactly matches $\mu_i$. If the moments $m_{\alpha_i}(\mathbf K_i)$ of the marginal distribution $\mu_i$ over $\mathbf K_i$ are either available analytically or efficiently computable for all $\alpha_i =1,\ldots, 2r$, then the above condition simplifies to the linear equations
\begin{equation}
	\label{eq:marginal-matching}
	\sum_{\substack{\alpha\in N(m,2r) \\ \alpha_i=0}} h_{\alpha} \prod_{j\neq i} m_{\alpha_j}(\mathbf{K}_j)=1\quad \text{and}
	\sum_{\substack{\alpha\in N(m,2r) \\ \alpha_i=\ell}} h_{\alpha} \prod_{j\neq i} m_{\alpha_j}(\mathbf{K}_j)=0 ~\forall \ell=1,\ldots, 2r, \ \forall i = 1,\ldots,m.
\end{equation}
Situations where the marginals of groups of random variables are known can be handled analogously. Note that when all marginals are known, there is only ambiguity about the dependence structure or {\em copula} of the components of $z$ \cite{Sklar1959}. Quantifying the worst-case copula amounts to solving a so-called Fr\'echet problem. In distributionally robust optimization, Fr\'echet problems with discrete marginals or approximate marginal matching conditions have been studied in \cite{Doan2012,Doan2015,VanParys2016}. \\[1ex]

Besides the ambiguity set $\mathcal P$, the constraint function $f$ also admits some generalizations that preserve the SDP-representability of the worst-case expectation in~\eqref{fminkreform2}.

\paragraph{Uncertainty quantification problems} If the constraint function $f$ in~\eqref{fminkreform2} is given by $f(x,z)=\mathbf 1_{z\in\mathbf C}$ for some Borel set $\mathbf C\subset \mathbb R^m$, then the worst-case expectation reduces to the worst-case probability of the set $\mathbf C$. Moreover, if the moments $m_\alpha(\mathbf K\cap \mathbf C)$ of the reference measure $\mu$ over $\mathbf K\cap \mathbf C$ are either available analytically or efficiently computable for all $\alpha \in N(m,2r)$, then the worst-case probability can be computed by solving a variant of the SDP~\eqref{eqSDP} with the alternative objective function
\[
	\sum_{\alpha\in N(m,2r)}h_{\alpha}m_{\alpha}(\mathbf{K}\cap\mathbf C).
\]

\section{Numerical experiments} \label{sec:results}
In the following we will exemplify the proposed approach to distributionally robust optimization in the context of financial portfolio analysis (Section~\ref{sec:portfolio}) and risk aggregation (Section~\ref{sec:Frechet}).

\subsection{Portfolio analysis} \label{sec:portfolio}
Consider a portfolio optimization problem, where the decision vector $x \in \mathbb{R}^n$ captures the percentage weights of the initial capital allocated to $n$ different assets. By definition, one thus has $x_i \in [0,1]$ for all $i= 1,\ldots,n$ and $\sum_{i} x_i = 1$. We assume that the asset returns $r_i = (u_i + l_i) / 2 + z_i (u_i - l_i) / 2$ depend linearly on some uncertain risk factors $z_i \in [-1,1]$ for all $i=1,\ldots, n$, where $l_i$ and $u_i$ represent known upper and lower bounds on the $i$-th return, respectively. In this framework, we denote by $z\in\mathbb R^n$ the vector of all risk factors and by $\mathbf K=[-1,1]^n$ its support. Moreover, the portfolio return can be expressed as
$$
f(x,z) = \sum\limits_{i=1}^n x_i \cdot ((u_i + l_i) / 2 + z_i (u_i - l_i) / 2).
$$
Unless otherwise stated, we set $\mu$ to the Lebesgue measure on $\mathbb R^n$. Modeling the probability density functions as SOS polynomials allows to account for various statistical properties and stylized facts of real asset returns as described in \cite{Cont2001}. For example, the proposed approach can conveniently capture \emph{gain\ loss asymmetry}, {\em i.e.}, the observation that large drawdowns in stock prices and stock index values are more common than equally large upward movements. This feature can be modeled by assigning a higher probability to an individual asset's large upward returns than to its low downward returns. Specifically, the ambiguity set may include the conditions $\mathbb{P}_{z \sim (\mathbf{K},h)}(z_i \leq a_i ) = \gamma_1$ and $\mathbb{P}_{z \sim (\mathbf{K},h)}(z_i \geq b_i ) = \gamma_2$ for some thresholds $a_i<b_i$ and confidence levels $\gamma_1 > \gamma_2$. %, where
%\begin{align*}
%\mathbb{P}(z_i \le a_i ) & = \sum_{\alpha\in N(m,2r)} h_{\alpha}m_{\alpha}(\mathbf{K} \cap \{z: \ z_i \leq a_i \}), \\
% \mathbb{P}(z_i \ge b_i ) & = \sum_{\alpha\in N(m,2r)} h_{\alpha}m_{\alpha}(\mathbf{K} \cap \{z: \ z_i \geq b_i \}).
%\end{align*}

Similarly, our approach can handle \emph{correlations of extreme returns}. As pointed in \cite{Cont2001}, in spite of the widespread use of the covariance matrix, `in circumstances when stock prices undergo large fluctuations [...], a more relevant quantity is the conditional probability of a large (negative) return in one stock given a large negative movement in another stock.' An example constraint on the conditional probability of one asset's low performance given another assets' lower performance is $\mathbb{P}_{z \sim (\mathbf{K},h)}(z_i \leq \underline{r}_i| z_j \leq \underline{r}_j ) \leq \gamma$, where $\underline r_i$ and $\underline r_j$ are given thresholds, while $\gamma$ is a confidence level.

In this numerical experiment we evaluate the probability that the return of a fixed portfolio $x$ materializes below a prescribed threshold $\underline r$, that is, we evaluate the worst case of the probability
\[
	\mathbb{P}_{z \sim (\mathbf{K},h)} \left( r(x,z) \leq \underline{r} \right)
\]
over an ambiguity set $\mathcal P$ of the form~\eqref{q:ambiguity-set} with the additional moment constraints ${\mathbb E}_{z \sim (\mathbf{K},h)} [z^{\beta_i}]=\gamma_i$ for some given exponents $\beta_i\in \mathbb N_0^n$ and targets $\gamma_i\in\mathbb R$ for $i=1,\ldots,p$. This corresponds to computing the integral of the density function over the knapsack polytope $\mathbf{K} \cap \mathbf{A}(x,u,l, \underline{r} )$, where
$$
\mathbf{A}(x,u,l, \underline{r} ) = \left\{ z\in\mathbb R^n: \ \sum_{i=1}^n x_i (u_i - l_i) z_i /2 \leq \underline{r} - \sum_{i=1}^n x_i (u_i + l_i) / 2  \right\}
$$
represents a halfspace in $\mathbb R^n$ that depends on the fixed portfolio $x$, the return bounds $l=(l_1,\ldots l_n)$ and $u=(u_1,\ldots,u_n)$, and the threshold $\underline r$. To formulate this problem as an SDP, we first need to compute the moments of the monomials with respect to the Lebesgue measure over the given knapsack polytope by using the results of \cite{Marichal2008}. The worst-case probability problem can then be reformulated as the SDP
\begin{align}
\label{eq:sdp2}
\begin{array}{cl}
\displaystyle\sup\limits_{h(z)} \ & \displaystyle\sum_{\alpha\in N(n,2r)} h_{\alpha}m_{\alpha}(\mathbf{K} \cap \mathbf{A}(x,u,l, \underline{r} )) \\[2ex]
\text{s.t.} \  & \displaystyle\sum_{\alpha\in N(n,2r)}h_{\alpha}m_{\alpha}(\mathbf{K})=1, \\[2ex]
& \displaystyle\sum_{\alpha\in N(n,2r)}h_{\alpha}m_{\alpha+\beta_i}(\mathbf{K})=\gamma_i \quad \forall i = 1,\ldots,p, \\[2ex]
& \displaystyle\sum_{\alpha\in N(n,2r)}h_{\alpha}z^{\alpha}\in\Sigma[z]_r.
\end{array}
\end{align}
In the numerical experiment we assume that there are $n = 2$ assets with lower and upper return bounds $l = (0.8, 0.7)^\top$ and $u = (1.2, 1.3)^\top$, respectively. We evaluate the probability that the return of the fixed portfolio $x = (0.75, 0.25)^\top$ falls below the threshold $\underline{r} = 0.9$ (the minimum possible return of the portfolio is 0.775). We assume that the only known moment information about the asset returns is that their means both vanish, that is, we set $p=2$, $\beta_1=(1,0)$, $\beta_2=(0.1)$ and $\gamma_1=\gamma_2=0$. Table~\ref{tab.portfolio.results.heuristic} reports the {\em exact} optimal values of the SDP~\eqref{eq:sdp2} for $r = 1,\ldots,12$ ($R=1$). The value in the last row of the table (labeled $r=\infty$) provides the worst-case probability across {\em all} distributions satisfying the prescribed moment conditions (not only those with a polynomial density) and was computed using the methods described in \cite{Hanasusanto2015}. In this case, one can also show that there exists a worst-case distribution with only two atoms. It assigns probability $0.31$ to the scenario $z=(1,1)^\top$ and probability $0.69$ to the scenario $z=(0.28,0.28)^\top$. We further computed the worst-case probabilities {\em approximately} by using Algorithm~\ref{alg:heuristic} from Section~\ref{section.heuristic} for SOS parameters $r = 1,\ldots,6$ and for up to $R=14$ iterations. The results are given in Table~\ref{tab.portfolio.results.heuristic} ($R>1$).

\begin{table}
\centering
\caption{Worst-case probabilities for the portfolio return falling below $\underline r$ computed by directly solving the SDP~\eqref{eq:sdp2} ($R=1$) and by using Algorithm~\ref{alg:heuristic} ($R>1$). Missing values indicate the occurrence of numerical instability.}
\label{tab.portfolio.results.heuristic}
\begin{tabular}{|c|cccccccccccccc|} \hline
%\multirow{2}{*}{$r$}
& \multicolumn{14}{c|}{$R$} \\
   $r$ & 1 & 2 & 3 & 4 & 5 & 6 & 7 & 8 & 9 & 10 & 11 & 12 & 13 & 14 \\ \hline
   0 & 0.17 &&&&&&&&&&&&& \\
   1 & 0.39 & 0.47 & 0.49 & 0.50 & 0.50 & 0.50 & 0.50 & 0.50 & 0.50 & 0.50 & 0.50 & 0.50 & 0.50 & 0.57 \\
   2 & 0.48 & 0.52 & 0.54 & 0.55 & 0.56 & 0.57 & 0.60 &-&-&-&-&-&-&-\\
   3 & 0.50 & 0.56 & 0.59 & 0.62  &-&-&-&-&-&-&-&-&-&- \\
   4 & 0.53 & 0.58 & 0.61  &-&-&-&-&-&-&-&-&-&-&-\\
   5 & 0.55 & 0.59 &-&-&-&-&-&-&-&-&-&-&-&-\\
   6 & 0.56 & 0.60 &-&-&-&-&-&-&-&-&-&-&-&-\\
   7 & 0.58 &-&-&-&-&-&-&-&-&-&-&-&-&-\\
   8 & 0.59 &-&-&-&-&-&-&-&-&-&-&-&-&-\\
   9 & 0.59 &-&-&-&-&-&-&-&-&-&-&-&-&-\\
  10 & 0.60 &-&-&-&-&-&-&-&-&-&-&-&-&-\\
  11 & 0.61 &-&-&-&-&-&-&-&-&-&-&-&-&-\\
  12 & 0.61 &-&-&-&-&-&-&-&-&-&-&-&-&-\\ \hline
  $\infty$ & 0.69 &-&-&-&-&-&-&-&-&-&-&-&-&-\\ \hline
\end{tabular}
%\label{tab.portfolio.results.heuristic}
%\begin{tabular}{cc|cccccccc} \hline
%\multicolumn{2}{c|}{Direct solution} & \multicolumn{7}{c}{Algorithm~\ref{alg:heuristic}} \\ \hline
%$r$ & & $R$: $\downarrow$ / $r$: $\rightarrow$ & 1 & 2 & 3 & 4 & 5 & 6 \\ \hline
%   0 & 0.1736 & 1 & 0.3943 & 0.4822 & 0.4988 & 0.5252 & 0.5549 & 0.5555 \\
%   1 & 0.3943 & 2 & 0.4712 & 0.5223 & 0.5626 & 0.5794 & 0.5933 & 0.5969 \\
%   2 & 0.4824 & 3 & 0.4927 & 0.5392 & 0.5864 & 0.6079 & \\
%   3 & 0.4988 & 4 & 0.4983 & 0.5516 & 0.6212   \\
%   4 & 0.5249 & 5 & 0.4996 & 0.5624  \\
%   5 & 0.5419 & 6 & 0.4999 & 0.5719  \\
%   6 & 0.5641 & 7 & 0.5000 & 0.5997  \\
%   7 & 0.5755 & 8 & 0.5000  \\
%   8 & 0.5889 & 9 & 0.5000  \\
%   9 & 0.5947 & 10 & 0.5000   \\
%  10 & 0.6023 & 11 & 0.5000   \\
%  11 & 0.6090 & 12 & 0.5008  \\
%  12 & 0.6142 & 12 & 0.5008  \\
%  13 & & 13 & 0.5019 \\
%  14 & & 14 & 0.5730 \\ \hline
%  $\infty$ & 0.6923  \\ \hline
%\end{tabular}
\end{table}

\subsection{Risk aggregation} \label{sec:Frechet}

In the second experiment we study the risk aggregation problem of an insurer holding a portfolio of different random losses $z_i$, $i=1,\ldots, n$, corresponding to different types of insurance claims, {\em e.g.}, life, vehicle, health or home insurance policies, etc. Inspired by \cite[\S~6]{VanParys2016}, we aim to estimate the worst-case probability that the sum of the $n$ losses exceeds a critical threshold $b=10$ beyond which the insurance company would be driven into illiquidity. Formally, we aim to maximize
\begin{equation}\label{eg:worstcase prob losses}
\mathbb{P}_{z \sim (\mathbf{K},h)} \left( z_1 + \ldots + z_n \geq b \right)
\end{equation}
across all distributions in an ambiguity set $\mathcal P$, which reflects the prior distributional information available to the insurer. We will consider different models for the domain $\bK$ of $z=(z_1,\ldots, z_n)$, the reference measure $\mu$ on $\bK$ and the ambiguity set $\mathcal P$. Throughout the experiments we will always assume that the reference measure is separable with respect to the losses, that is, we assume that
\begin{equation*}
d\mu(z) = \varrho_1(z_1) \cdots \varrho_n(z_n) dz,
\end{equation*}
where $\varrho_i$ denotes a given density function (with respect to the Lebesgue measure) of the random variables $z_i$ for each $i=1,\ldots, n$. We will consider the following complementary settings:

\begin{enumerate}
\item {\em Lognormal densities:}
We set $\bK = \mathbb{R}^n_+$ and let $\varrho_i$ be a lognormal density function defined earlier in \eqref{eq:lognormal0}, but repeated here for convenience:
\begin{equation}
	\label{eq:lognormal}
	\varrho_i(z_i) = \frac{1}{z_i v_i \sqrt{2\pi}} \exp\left( - \frac{(\log(z_i) - \bar z_i)^2}{2 v_i^2} \right),
\end{equation}
where $\bar z_i$ and $v_i$ represent prescribed location and scale parameters, $i = 1,\ldots,n$.
\item {\em Exponential densities:}
We set $\bK = \mathbb{R}^n_+$ and let $\varrho_i$ be the exponential density function with unit rate parameter defined through $\varrho_i(z_i)=\exp(-z_i)$, $i=1,\ldots, n$. The resulting reference measure is intimately related to the orthogonal Laguerre polynomials.
\item {\em Uniform densities:}
We set $\bK = [0,M]^n$ for some constant $M>0$ and let $\varrho_i$ be the uniform density function defined through $\varrho_i(z_i)=1/M$, $i=1,\ldots, n$. Note that under this choice the reference measure is proportional to the Lebesgue measure.
\end{enumerate}

In order to reformulate the risk aggregation problem as a tractable SDP, we need the moments of the reference measure $\mu$ over the hypercube $\mathbf K$ and over the knapsack polytope $\mathbf K\cap\mathbf C$, where
\[
	\mathbf C = \{z \in \mathbb{R}^n: \ \ z_1 + \ldots + z_n \geq b \}.
\]
For all classes of density functions described above, the moments of $\mu$ are indeed accessible. Specifically, under the lognormal densities, the moments of $\mu$ over $\mathbf K$ are
given by \eqref{eq:lognormal moments}, and repeated here for convenience:
$$
m_{\alpha}(\mathbf K) = \int_{\mathbf K} \left( \prod_{i=1}^n z_i^{\alpha_i} \right) \prod_{i=1}^n \frac{1}{z_i v_i \sqrt{2\pi}}
 \exp \left( - \frac{(\log(z_i) - \bar z_i)^2}{2 v_i^2} \right) dz = \prod_{i=1}^n \exp(\alpha_i \bar z_i + (\alpha_i v_i)^2 / 2).
$$
Moreover, the moments of $\mu$ over $\mathbf K\cap\mathbf C$ can be expressed as
\begin{align*}
m_{\alpha}(\mathbf K\cap \mathbf C) & =m_{\alpha}(\mathbf K) -m_{\alpha}(\mathbf K\backslash \mathbf C)\\
& =m_{\alpha}(\mathbf K) - \int_{\mathbf K\backslash \mathbf C} \left( \prod_{i=1}^n z_i^{\alpha_i} \right) \prod_{i=1}^n \frac{1}{z_i v_i \sqrt{2\pi}} \exp \left( - \frac{(\log(z_i) - \bar z_i)^2}{2 v_i^2} \right) dz.
\end{align*}
To evaluate the integral in the last expression, we use the  MATLAB routine {\tt adsimp($\cdot$)} from \cite{Genz2003}. Furthermore, under the exponential and the uniform densities, the moments of the reference measure $\mu$ over $\mathbf K$ and $\mathbf K\cap\mathbf C$ are all available in closed form. %Likewise, under the uniform densitiesall the required moments are available in closed form.

\begin{figure}[h!]
  \input{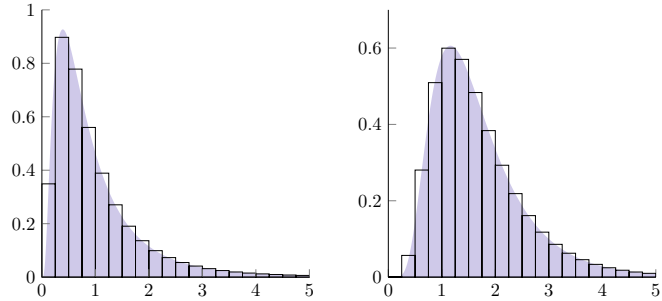}
\caption{Histograms of the lognormal marginal distributions of $z_1$ (left) and $z_2$ (right). \label{fig:histograms}}
\end{figure}

We assume that the insurance company is able to estimate the marginal distributions of the individual losses either exactly or approximately by using a combination of statistical analysis and probabilistic modeling. However, the insurer has no information about the underlying copula. This type of distributional information is often justified in practice because obtaining reliable marginal information requires significantly less data than obtaining exact dependence structures; see, {\em e.g.},~\cite{McNeil2015}. Throughout the experiment we assume that there are $n=2$ random losses governed by lognormal probability density functions of the form~\eqref{eq:lognormal} with parameters $\bar z_1 = -0.3$,  $\bar z_2 = 0.4$, $v_1 = 0.8$ and $v_2 = 0.5$. The ambiguity set $\mathcal P$ then contains all distributions of the form $h\cdot d\mu$, $h\in \Sigma[z]_r$, under which the marginals of the losses follow the prescribed lognormal distributions either exactly or approximately. More precisely, we model the marginal distributional information as follows:

%Following \cite{VanParys2016}, we take these desirable properties to mean that the marginals of the joint distribution that we obtain should be close to the functions $f_i$, i.e.\ we view  $f_i$ as the `true' marginal distribution of $z_i$. We can model the `close to the required marginal' requirement in several ways:
 \begin{enumerate}
 \item {\em Marginal distribution matching:} The lognormal distributions of the individual losses are matched {\em exactly} by any distribution $h\cdot d\mu$ in the ambiguity set. This can be achieved by defining the reference measure $\mu$ as
 the product of the marginal lognormal distributions and by requiring that $h$ satisfies \eqref{eq:marginal-matching}. Note that under the alternative reference measures corresponding to the exponential or uniform density functions, lognormal marginals cannot be matched exactly with polynomial densities of any degrees. Note also that an exact matching of (non-discrete) marginal distributions cannot be enforced with the existing numerical techniques for solving Fr\'echet problems proposed in \cite{Doan2012,Doan2015,VanParys2016}.
 \item {\em Marginal moment matching:}
The marginals of the individual losses have the same moments of order 0, 1 or 2 as the prescribed lognormal distributions. Note that this kind of moment matching can be enforced under any of the reference measures corresponding to lognormal, exponential or uniform density functions. Moreover, moment matching is also catered for in \cite{VanParys2016} bar the extra requirement that the joint distribution of the losses must have an SOS polynomial density.
 \item {\em Marginal histogram matching:}
 We may associate a histogram with each marginal lognormal distribution as illustrated in Figure~\ref{fig:histograms} and require that the marginals of the losses under the joint distribution $h\cdot d\mu$ have the same histograms. This condition can be enforced under any of the reference measures corresponding to lognormal, exponential or uniform density functions. In the numerical experiments, we use histograms with 20 bins of width 0.25 starting at the origin. Histogram matching is also envisaged in~\cite{VanParys2016}.
  \end{enumerate}

For $\bK=\mathbb R_+^n$ and the reference measure corresponding to the lognormal density functions, the worst-case values of the probability~\eqref{eg:worstcase prob losses} are reported in Table~\ref{tab.3}. Results are shown for $r\leq 5$, which corresponds to polynomial densities of degrees at most $10$. The last row of the table ($r=\infty$) provides the worst-case probabilities across {\em all} distributions satisfying the prescribed moment or histogram conditions (not only those with a polynomial density) and was computed using the methods described in \cite{VanParys2016}. Note that under moment matching up to order $2$, the worst-case probability for $r=5$ amounts to $0.0021$, as opposed to the much higher probability of $0.0615$ obtained with the approach from \cite{VanParys2016}. A similar observation holds for histogram matching. The requirement that the distributions in the ambiguity set be sufficiently regular in the sense that they admit a polynomial density function with respect to the reference measure is therefore restrictive and effectively rules out pathological discrete worst-case distributions. Moreover, the worst-case probabilities under exact distribution matching and under histogram matching are of the same order of magnitude for all $r\le 5$ but significantly smaller than the worst-case pobability under histogram matching for $r=\infty$. A key question to be asked in practice is thus whether one deems the class of distributions $h\cdot d\mu$ with $h \in \Sigma[z]_r$ to be rich enough to contain all `reasonable' distributions.

\begin{table}[h!]
\centering
\caption{Worst-case probabilities for the lognormal reference measure.} \label{tab.3}
\begin{tabular}{|c|ccc|c|c|} \hline
  & \multicolumn{3}{c|}{Moment matching up to order} & Histogram & Distribution\\
$r$  & 0 & 1 & 2 & matching &  matching \\ \hline
0 & 0.0017 & 0.0017 & 0.0017 & 0.0017 & 0.0017 \\
1 & 0.1432 & 0.0042 & 0.0017 & 0.0017 & 0.0017 \\
2 & 0.8255 & 0.0106 & 0.0020 & 0.0019 & 0.0018 \\
3 & 0.9982 & 0.0114 & 0.0021 & 0.0022 & 0.0019 \\
4 & 1.0000 & 0.0117 & 0.0021 & 0.0026 & 0.0023 \\
5 & 1.0000 & 0.0118 & 0.0021 & 0.0026 & 0.0023 \\ \hline
$\infty$ & 1.0000 & 1.0000 & 0.0615 & 0.0198 & n/a \\ \hline
\end{tabular}
\end{table}

\begin{table}[h!]
\centering
\caption{Worst-case probabilities for the exponential reference measure.} \label{tab.Laguerre}
\begin{tabular}{|c|ccc|ccc|} \hline
 & \multicolumn{3}{c|}{Moment matching up to order} & \multicolumn{3}{c|}{Histogram matching} \\
$r$  & 0 & 1 & 2 &  $\ell_1$-dist. $\leq 0.1$ & $\ell_1$-dist. $\leq 0.05$ & $\ell_1$-dist. $\leq 0.02$\\ \hline
  0 & 0.0005 &      - &      -  &      - &      -  &      - \\
  1 & 0.0214 & 0.0147 &      -  &      - &      -  &      - \\
  2 & 0.2058 & 0.0823 &      -  &      - &      -  &      - \\
  3 & 0.6481 & 0.1484 &      -  &      - &      -  &      - \\
  4 & 0.9393 & 0.1497 & 0.0086  &      - &      -  &      - \\
  5 & 0.9953 & 0.1699 & 0.0104  &      - &      -  &      - \\
  6 & 0.9998 & 0.1709 & 0.0139  &      - &      -  &      - \\
  7 & 1.0000 & 0.1800 & 0.0158  &      - &      -  &      - \\
  8 & 1.0000 & 0.1860 & 0.0182  & 0.0802 &      -  &      - \\
  9 & 1.0000 & 0.1862 & 0.0207  & 0.1076 &      -  &      - \\
 10 & 1.0000 & 0.1928 & 0.0224  & 0.1144 & 0.0515  &      - \\
 11 & 1.0000 & 0.1968 & 0.0244  & 0.1156 & 0.0633  & 0.0204 \\
 12 & 1.0000 & 0.1971 & 0.0262  & 0.1160 & 0.0652  & 0.0320 \\ \hline
  $\infty$ & 1.0000 & 1.0000 & 0.0615 & n/a & n/a & n/a \\ \hline
% 13 & 1.888838 & 0.318810 & 0.030506 & 0.006182 & 0.010915 \\ \hline
\end{tabular}
\end{table}

%\begin{table}[h!]
%\centering
%\caption{Worst-case probabilities for exponential reference measure.} \label{tab.Laguerre}
%\begin{tabular}{|c|ccc|cc|} \hline
% & \multicolumn{3}{c|}{Moment matching up to order} & \multicolumn{2}{c|}{Histogram matching}  \\
%$r$  & 0 & 1 & 2  & Prob. & $\ell_1$-dist. \\ \hline
%  0 & 0.0005 &        - &        - & 0.0005 & 1.2232 \\
%  1 & 0.0214 & 0.0147 &        - & 0.0049 & 0.8696 \\
%  2 & 0.2058 & 0.0823 &        - & 0.0491 & 0.4906 \\
%  3 & 0.6481 & 0.1484 &        - & 0.0378 & 0.3517 \\
%  4 & 0.9393 & 0.1497 & 0.0086 & 0.0395 & 0.2812 \\
%  5 & 0.9953 & 0.1699 & 0.0104 & 0.0391 & 0.2198 \\
%  6 & 0.9998 & 0.1709 & 0.0139 & 0.0428 & 0.1725 \\
%  7 & 1.0000 & 0.1800 & 0.0158 & 0.0460 & 0.1214 \\
%  8 & 1.0000 & 0.1860 & 0.0182 & 0.0349 & 0.0849 \\
%  9 & 1.0000 & 0.1862 & 0.0207 & 0.0274 & 0.0606 \\
% 10 & 1.0000 & 0.1928 & 0.0224 & 0.0267 & 0.0422 \\
% 11 & 1.0000 & 0.1968 & 0.0244 & 0.0129 & 0.0275 \\
% 12 & 1.0000 & 0.1971 & 0.0262 & 0.0078 & 0.0177 \\ \hline
%  $\infty$ & 1.0000 & 1.0000 & 0.0615 & 0.0198 & 0.0000 \\ \hline
%% 13 & 1.888838 & 0.318810 & 0.030506 & 0.006182 & 0.010915 \\ \hline
%\end{tabular}
%\end{table}

Table~\ref{tab.Laguerre} reports the worst-case probabilities corresponding to the reference measure on $\bK=\mathbb R_+^n$ induced by the exponential density functions. For low values of $r$, the polynomial densities lack the necessary degrees of freedom to match all imposed moment constraints. In these situations, the worst-case probability problem becomes infeasible (indicated by `-'). When feasible, however, we managed to solve the problem for $r$ up to 12. The density functions corresponding to large values of $r$ are highly flexible and thus result in worst-case probabilities that are closer to those obtained by the benchmark method from~\cite{VanParys2016}, which relaxes the restriction to a subspace of polynomial densities. Similar phenomena are also observed in the context of histogram matching. It was impossible to match the prescribed histogram probabilities exactly for all $r\le 12$. We thus relaxed the histogram matching conditions in the definition of the ambiguity set to allow for densities whose implied marginal histograms are within a prescribed $\ell_1$-distance from the target histograms. This approximate histogram matching condition is easily captured in our framework and gives rise to a few extra linear constraints on the coefficients of the polynomial density function. Table~\ref{tab.Laguerre} reports the worst-case probabilities for three different tolerances on the histogram mismatch in terms of the $\ell_1$-distance. We observe that the resulting worst-case probabilities are significantly larger than those obtained under the lognormal reference measure and increase with the $\ell_1$-tolerance. %This gives an indication of how the worst-case probability changes relative to inaccuracies in the histogram.

Finally, Table~\ref{tab.Lebesgue} reports the worst-case probabilities corresponding to the uniform reference measure on $\bK = [0,10]^2$. The results are qualitatively similar to those of Table~\ref{tab.Laguerre}, but they also show that the choice of the reference measure plays an important role when $r$ is small.
\begin{table}[h!]
\centering
\caption{Worst-case probabilities for the uniform reference measure.} \label{tab.Lebesgue}
\begin{tabular}{|c|ccc|ccc|} \hline
 & \multicolumn{3}{c|}{Moment matching up to order} &  \multicolumn{3}{c|}{Histogram matching} \\
$r$  & 0 & 1 & 2 & $\ell_1$-dist. $\leq 0.1$ & $\ell_1$-dist. $\leq 0.05$ & $\ell_1$-dist. $\leq 0.02$ \\ \hline
0 & 0.5000 & - & - & -      & -      & -  \\
1 & 0.9082 & - & - & -      & -      & -  \\
2 & 0.9933 & - & - & -      & -      & -  \\
3 & 0.9997 & 0.0304 & - & -      & -      & -  \\
4 & 1.0000 & 0.1035 & - & -      & -      & -  \\
5 & 1.0000 & 0.1340 & - & -      & -      & -  \\
6 & 1.0000 & 0.1612 & 0.0089 & -      & -      & -  \\
  7 & 1.0000 & 0.1783 & 0.0166 & -      & -      & -  \\
  8 & 1.0000 & 0.1935 & 0.0192 & -      & -      & - \\
  9 & 1.0000 & 0.2042 & 0.0216 & 0.0738 & -      & - \\
 10 & 1.0000 & 0.2133 & 0.0274 & 0.1066 & 0.0407 & -  \\
 11 & 1.0000 & 0.2202 & 0.0292 & 0.1142 & 0.0609 & -  \\
 12 & 1.0000 & 0.2274 & 0.0311 & 0.1163 & 0.0653 & 0.0178 \\ \hline
 $\infty$ & 1.0000 & 1.0000 & 0.0615 &  n/a & n/a & n/a \\ \hline
\end{tabular}
\end{table}

%\begin{table}[h!]
%\centering
%\caption{Worst-case probabilities for uniform reference measure.} \label{tab.Lebesgue}
%\begin{tabular}{|c|ccc|cc|} \hline
% & \multicolumn{3}{c|}{Moment matching up to order} & \multicolumn{2}{c|}{Histogram matching} \\
%$r$  & 0 & 1 & 2  & Prob. & $\ell_1$-dist.\\ \hline
%0 & 0.5000 & - & - & 0.5000 & 1.7107 \\
%1 & 0.9082 & - & - & 0.3696 & 1.3029 \\
%2 & 0.9933 & - & - & 0.1557 & 1.0081 \\
%3 & 0.9997 & 0.0304 & - & 0.1554 & 0.8494 \\
%4 & 1.0000 & 0.1035 & - & 0.0915 & 0.6758 \\
%5 & 1.0000 & 0.1340 & - & 0.0849 & 0.4439 \\
%6 & 1.0000 & 0.1612 & 0.0089 & 0.0473 & 0.2841 \\
%  7 & 1.0000 & 0.1783 & 0.0166 & 0.0268 & 0.2048 \\
%  8 & 1.0000 & 0.1935 & 0.0192 & 0.0283 & 0.1446 \\
%  9 & 1.0000 & 0.2042 & 0.0216 & 0.0232 & 0.0809 \\
% 10 & 1.0000 & 0.2133 & 0.0274 & 0.0131 & 0.0410 \\
% 11 & 1.0000 & 0.2202 & 0.0292 & 0.0079 & 0.0173 \\
% 12 & 1.0000 & 0.2274 & 0.0311 & 0.0041 & 0.0041 \\ \hline
% $\infty$ & 1.0000 & 1.0000 & 0.0616 & 0.0199 &0.0000 \\ \hline
%\end{tabular}
%\end{table}

\section{Conclusions}
\label{sec:conclusions}
In this paper, we present first steps towards using SOS polynomial densities in distributionally robust optimization for problems that display a polynomial dependence on the uncertain parameters. The main advantages of this approach may be summarized as follows:
\begin{enumerate}
\item
The proposed framework is tractable (in the sense of polynomial-time
solvability) for SOS density functions of any fixed degree.
\item
The approach offers considerable modeling flexibility. Specifically, one may conveniently encode various salient features of the unknown distribution of the uncertain parameters trough linear constraints and/or linear matrix inequalities.
\item
In the limit as the degree of the SOS density functions tends to infinity, one recovers the usual robust counterpart or generalized moment problem. One may therefore view the degree of the density as a tuning parameter that captures the model's  `level of conservativeness.'
\end{enumerate}

The approach also suffers from shortcomings that necessitate further work and insights:
\begin{enumerate}
\item
The approach is not applicable to objective or constraint functions that display a general (decision-dependent) piecewise polynomial dependence on the uncertain parameters as is the case for the recourse functions of linear two-stage stochastic programs.
%In applications, one will have to give statistical justification for modeling the distribution of the uncertain parameters by a polynomial density. In particular, one has to argue that this class is rich enough for the specific application.
\item
The proposed distributionally robust optimization problems can be reduced to generalized eigenvalue problems or even semidefinite programs of large sizes that are often poorly conditioned. We have introduced a heuristic solution procedure in Section~\ref{section.heuristic} as a practical remedy, but additional work is required to make the approach more scalable.
\end{enumerate}

\paragraph{\bf Acknowledgements}
Etienne de Klerk would like to thank Dorota Kurowicka and Jean Bernard Lasserre for valuable discussions and references. Daniel Kuhn gratefully acknowledges financial support from the Swiss National Science Foundation under grant BSCGI0\_157733.


\begin{thebibliography}{10}

\bibitem{BenTal2009}
Ben-Tal, A., {El Ghaoui}, L., Nemirovski, A., \emph{Robust Optimization}, Princeton University Press (2009).

\bibitem{BenTal2013}
Ben-Tal, A. {den Hertog}, D., {de Waegenaere}, A., Melenberg, B., Rennen, G. Robust solutions of optimization problems affected by uncertain probabilities. \emph{Management Science} 59(2), 341--357 (2013).

\bibitem{Bertsimas2002}
Bertsimas, D., Popescu, I. On the relation between option and stock prices: a convex optimization approach, \emph{Operations Research} 50(2), 358--374 (2002).

\bibitem{Bertsimas2005}
Bertsimas, D., Popescu, I. \emph{Optimal inequalities in probability theory: a convex optimization approach}, \emph{SIAM Journal on Optimization} 15(3), 780--804 (2005).

\bibitem{Birge1997}
Birge, J.R., Louveaux, F. \emph{Introduction to Stochastic Programming}, Springer (1997).

\bibitem{Cont2001}
Cont, R. Empirical properties of asset returns: stylized facts and statistical issues. \emph{Quantitative
Finance} 1, 223--236 (2001).

\bibitem{KL_MOR_2017}
{de Klerk}, E., Laurent, M. Comparison of Lasserre's measure-based  bounds for polynomial optimization to bounds obtained by simulated annealing. \emph{Mathematics of Operations Research}, to appear. Preprint available at \url{http://arxiv.org/abs/1703.00744}

\bibitem{KL_2018}
{de Klerk}, E., Laurent, M. Worst-case examples for Lasserre's measure--based hierarchy for polynomial optimization on the hypercube (2018). Preprint available at \url{http://arxiv.org/abs/1804.05524}

\bibitem{KLS_MPA}
{de Klerk}, E., Laurent, M., Sun, Z. Convergence analysis for Lasserre's measure-based hierarchy of upper bounds for polynomial optimization. \textit{Mathematical Programming Series A} 162(1), 363--392 (2017).

\bibitem{Dantzig1955}
Dantzig, G.B. Linear programming under uncertainty. \emph{Management Science} 1(3-4), 197--206 (1955).

\bibitem{Delage2010}
Delage, E., Ye, Y. Distributionally robust optimization under moment uncertainty with application to data-driven problems. \emph{Operations Research} 58(3), 595--612 (2010).

\bibitem{Doan2015}
 Doan, X.V., Li, X., Natarajan, K. Robustness to dependency in portfolio optimization using overlapping marginals. \emph{Operations Research} 63(6), 1468--1488 (2015).

\bibitem{Doan2012}
Doan X.V., Natarajan K. On the complexity of nonoverlapping multivariate marginal bounds for probabilistic combinatorial optimization problems. \emph{Operations Research} 60(1),138--49 (2012).

\bibitem{Genz2003}
Genz, A., Cools, R. An adaptive numerical cubature algorithm for simplices. \emph{ACM Transactions on Mathematical Software} 29(3), 297--308 (2003).

\bibitem{Goh2010}
Goh, J., Sim, M. Distributionally robust optimization and its tractable approximations, \emph{Operationis Research} 58(4), 902--917 (2010).

\bibitem{GroetschelLovaszSchrijver1988a}
Gr\"otschel, M., Lov\'asz, L., Schrijver, A. {\em Geometric Algorithms and Combinatorial Optimization}. Springer (1988).

\bibitem{GM78}
Grundmann, A., Moeller, H.M. Invariant integration formulas for the $n$-simplex by combinatorial methods. {\em SIAM Journal on Numerical Analysis} 15, 282--290 (1978).

\bibitem{Hanasusanto2015c}
Hanasusanto, G.A., Roitch, V., Kuhn, D., Wiesemann, W. A distributionally robust perspective on uncertainty quantification and chance constrained programming. \emph{Mathematical Programming Series B} 151(1), 35--62 (2015).

\bibitem{Hanasusanto2015}
Hanasusanto, G.A., Roitch, V., Kuhn, D., Wiesemann, W. Ambiguous joint chance constraints under mean and dispersion information. \emph{Operations Research} 65(3), 751--767 (2017).

\bibitem{Hanasusanto2015b}
Hanasusanto, G.A., Kuhn, D., Wallace, S.W., Zymler, S. Distributionally robust multi-item newsvendor problems with multimodal demand distributions. \emph{Mathematical Programming Series A} 152(1), 1--32 (2015).

\bibitem{Kroo1999}
Kroo, A., Szil\'{a}rd, R. On Bernstein and Markov-type inequalities for multivariate polynomials on convex bodies. {\em Journal of Approximation Theory} 99(1), 134--152 (1999).

\bibitem{LZ01}
Lasserre, J.B., Zeron, E.S. Solving a class of multivariate integration problems via Laplace techniques. {\em Applicationes Mathematicae} 28(4), 391--405 (2001).

\bibitem{Las2008}
Lasserre, J.B. A semidefinite programming approach to the generalized problem of moments.
\emph{Mathematical Programming Series B} 112, 65--€"92 (2008).

\bibitem{Las11}
Lasserre, J.B. A new look at nonnegativity on closed sets and polynomial optimization. {\em SIAM Journal on Optimization} 21(3), 864--885 (2011).

\bibitem{Las12}
Lasserre, J.B.
The $\mathbf K$-moment problem for continuous linear functionals. \emph{Transactions of the American Mathematical Society} 365(5),  2489--2504 (2012).

\bibitem{Las2018}
Lasserre, J.B., Weisser, T. Representation of distributionally robust chance-constraints (2018). Preprint available at \url{http://arxiv.org/abs/1803.11500}

\bibitem{Li2017}
Li, B., Jiang, R., Mathieu, J.L. Ambiguous risk constraints with moment and unimodality information, \emph{Mathematical Programming Series A}, to appear (2017). Preprint available at \url{http://www.optimization-online.org/DB_FILE/2016/09/5635.pdf}

\bibitem{McNeil2015}
McNeil, A., Frey, R., Embrechts, P. \emph{Quantitative Risk Management: Concepts, Techniques and Tools}, Princeton University Press (2015).

\bibitem{Marichal2008}
Marichal,  J.-L.,  Mossinghof, M.J. Slices, slabs, and sections of the unit hypercube. \emph{ Online Journal of Analytic Combinatorics} 3, 1--11 (2008).

\bibitem{Mevissen2013}
Mevissen, M., Ragnoli, E., Yu, J.Y. Data-driven distributionally robust polynomial optimization, \emph{Advances in Neural Information Processing Systems (NIPS)} 26 (2013).

\bibitem{Esfahani2017}
{Mohajerin Esfahani}, P., Kuhn, D. Data-driven distributionally robust optimization using the Wasserstein metric: performance guarantees and tractable reformulations, \emph{Mathematical Programming Series A}, to appear (2017). Preprint available at \url{https://arxiv.org/abs/1505.05116}

\bibitem{Natarajan2009}
Natarajan, K., Pachamanova, D., Sim, M. Constructing risk measures from uncertainty sets, \emph{Operations Research} 57(5), 1129--1141 (2009).

\bibitem{Pflug2012}
Pflug, G.C., Pichler, A., Wozabal, D. The $1/N$ investment strategy is optimal under high model ambiguity, \emph{Journal of Banking \& Finance} 36(2), 410--417 (2012).

\bibitem{Pflug2007}
Pflug, G.C., Wozabal, D. Ambiguity in portfolio selection, \emph{Quantitative Finance} 7, 435--442 (2007).

\bibitem{Popescu2005}
Popescu, I. A semidefinite programming approach to optimal-moment bounds for convex classes of distributions, \emph{Mathematics of Operations Research} 30(3), 632--657 (2005).

\bibitem{Prekopa1995}
Pr\'{e}kopa, A. \emph{Stochastic Programming}, Kluwer Academic Publishers (1995).

\bibitem{Rockafellar}
Rockafellar, R.T. \emph{Convex Analysis}, Princeton University Press (1970).

\bibitem{rogosinski}
Rogosinski, W.W. Moments of non-negative mass, \textit{Proceedings of the Royal Society A} 245, 1--27 (1958).

\bibitem{Scarf1958}
Scarf, H. A min-max solution of an inventory problem, in H. Scarf, K. Arrow, and S. Karlin (Eds.), \emph{Studies in the Mathematical Theory of Inventory and Production}, Volume 10, pp.~201--209. Stanford University Press (1958).

\bibitem{Sklar1959}
Sklar, A. Fonctions de r\'{e}partition \`{a} $n$ dimensions et leurs marges, \emph{Publications de l'Institut de Statistique de L'Universit\'{e} de Paris} 8, 229--231 (1959).

\bibitem{Shapiro2001}
Shapiro, A. \emph{On duality theory of conic linear problems}, Semi-Infinite Programming: Recent Advances (M,{\'A}. Goberna and M.A. L{\'o}pez,
  eds.), Springer, 135--165 (2001).

\bibitem{Shapiro2009}
Shapiro, A., Dentcheva, D., Ruszczy\'nski, A. \emph{Lectures on Stochastic Programming: Modeling and Theory}, SIAM (2009).

\bibitem{VanParys2016}
{Van Parys}, B.P.G., Goulart, P.J., Embrechts, P. Fr\'{e}chet inequalities via convex optimization (2016). Preprint available at \url{http://www.optimization-online.org/DB_FILE/2016/07/5536.pdf}

\bibitem{VanParys2016b}
Van Parys, B.P.G., Goulart, P.J., Kuhn, D. Generalized Gauss inequalities via semidefinite programming. \emph{Mathematical Programming Series A} 156(1-2), 271--302 (2016).

\bibitem{Wiesemann2014}
Wiesemann, W., Kuhn, D., Sim, M. Distributionally robust convex optimization. \emph{Operations Research} 62(6), 1358--1376 (2014).

\bibitem{Zackova1966}
\v{Z}\'{a}\v{c}kov\'{a}, J. On minimax solutions of stochastic linear programming problems, \emph{\v Casopis pro p{\v{e}}stov\'an\'i matematiky} 91, 423--430 (1966).

\bibitem{Zuluaga2005}
Zuluaga, L., Pe\~na J.F. A conic programming approach to generalized Tchebycheff inequalities, \emph{Mathematics of Operations Research} 30(2), 369--388 (2005).

\end{thebibliography}
\end{document}